\documentclass[a4paper,english]{IEEEtran}
\usepackage[T1]{fontenc}
\usepackage[latin9]{inputenc}
\usepackage[active]{srcltx}
\usepackage{array}
\usepackage{multirow}
\usepackage{amsmath}
\usepackage{amssymb}
\usepackage{graphicx}
\usepackage{esint}

\usepackage[breaklinks]{hyperref}
\hypersetup{
	bookmarks=true,         % show bookmarks bar?
	colorlinks=true,       % false: boxed links; true: colored links
    linkcolor=blue,          % color of internal links 
    citecolor=green,        % color of links to bibliography
    filecolor=magenta,      % color of file links
    urlcolor=magenta           % color of external links
    }
\makeatletter

\providecommand{\tabularnewline}{\\}

\usepackage{mathtools}
\makeatother

\usepackage{babel}
\begin{document}
\global\long\def\Lone{L_{1}}
\global\long\def\Ltwo{L_{2}}
\global\long\def\v{\mathrm{v}}
\global\long\def\w{\mathrm{w}}
\global\long\def\f{f}
\global\long\def\part{T}
\global\long\def\x{\mathbf{x}}
\global\long\def\R{\mathbb{R}}
\global\long\def\vecspace{\mathrm{V}_{\part}}
\global\long\def\pip{\pi_{\part}}
\global\long\def\ortproj{P_{\part}}
\global\long\def\inner#1#2{\left\langle #1,#2\right\rangle }
\global\long\def\dif{\mathrm{d}}
\global\long\def\partOpt{\part^{\ast}}
\global\long\def\partQuasiOpt{\part^{(\ast)}}
\global\long\def\partUnif{\part_{U}}
\global\long\def\vecspaceOpt{\mathrm{V}_{\partOpt}}
\global\long\def\pipOpt{\pi_{\partOpt}}
\global\long\def\ortprojOpt{P_{\partOpt}}
\global\long\def\pipQuasiOpt{\pi_{\partQuasiOpt}}
\global\long\def\ortprojQuasiOpt{P_{\partQuasiOpt}}
\global\long\def\pipUnif{\pi_{\partUnif}}
\global\long\def\ortprojUnif{P_{\partUnif}}
\global\long\def\fsec{\f^{\prime\prime}}
\global\long\def\fimax{|\fsec_{i}|_{\max}}
\global\long\def\fZero{\fsec_{0x}}
\global\long\def\fOne{\fsec_{1x}}
\global\long\def\bb{\mathbf{b}}
\global\long\def\bg{\mathbf{g}}
\global\long\def\bv{\mathbf{v}}
\global\long\def\by{\mathbf{y}}
\global\long\def\mA{\mathtt{A}}
\global\long\def\mH{\mbox{H}}
\global\long\def\linei{\mathrm{line}_{i}}
\global\long\def\signum{\text{sign}}
\global\long\def\leftDeltaY#1{\Delta y_{#1}}
\global\long\def\rightDeltaY#1{\Delta y_{#1}}
\global\long\def\step{\mathbf{s}}
\global\long\def\fbest{\hat{\f}}
\global\long\def\cost{\text{cost}}
\global\long\def\vf{\mathbf{f}}

\title{Optimal polygonal $\Lone$ linearization and\\
fast interpolation of nonlinear systems}

\author{Guillermo~Gallego, Daniel~Berjón and~Narciso~García%
\thanks{This work has been partially supported by the Ministerio de Economía
y Competitividad of the Spanish Government under project TEC2010-20412
(Enhanced 3DTV). G. Gallego is supported by the Marie Curie - COFUND
Programme of the EU (Seventh Framework Programme).

G.~Gallego, D.~Berjón and N.~García are with Grupo de Tratamiento
de Imágenes (GTI), ETSI Telecomunicación, Universidad Politécnica
de Madrid, Madrid, Spain, e-mail: \{ggb,dbd,narciso\}@gti.ssr.upm.es.

%Copyright (c) 2014 IEEE. Personal use of this material is permitted.
%However, permission to use this material for any other purposes must
%be obtained from the IEEE by sending an email to pubs-permissions@ieee.org. %
%The final publication is available at IEEE via http://dx.doi.org/10.1109/TCSI.2014.2327313
% doi: 10.1109/TCSI.2014.2327313
Copyright (c) 2014 IEEE. Personal use of this material is permitted.
However, permission to use this material for any other purposes must
be obtained from the IEEE. 
\protect\href{http://dx.doi.org/10.1109/TCSI.2014.2327313}{DOI: 10.1109/TCSI.2014.2327313}
%The final publication is available at IEEE via http://dx.doi.org/10.1109/TCSI.2014.2327313
%
}}
\maketitle
\begin{abstract}
The analysis of complex nonlinear systems is often carried out using
simpler piecewise linear representations of them. A principled and
practical technique is proposed to linearize and evaluate arbitrary
continuous nonlinear functions using polygonal (continuous piecewise
linear) models under the L1 norm. A thorough error analysis is developed
to guide an optimal design of two kinds of polygonal approximations
in the asymptotic case of a large budget of evaluation subintervals
N. The method allows the user to obtain the level of linearization
(N) for a target approximation error and vice versa. It is suitable
for, but not limited to, an efficient implementation in modern Graphics
Processing Units (GPUs), allowing real-time performance of computationally
demanding applications. The quality and efficiency of the technique
has been measured in detail on two nonlinear functions that are widely
used in many areas of scientific computing and are expensive to evaluate.\end{abstract}
\begin{IEEEkeywords}
Piecewise linearization, numerical approximation and analysis, least-first-power,
optimization. 
\end{IEEEkeywords}

\section{Introduction}

\IEEEPARstart{T}{he} approximation of complex nonlinear systems
by simpler piecewise linear representations is a recurrent and attractive
task in many applications since the resulting simplified models have
lower complexity, fit into well established tools for linear systems
and are capable of representing arbitrary nonlinear mappings. Examples
include, among others, complexity reduction for finding the inverse
of nonlinear functions~\cite{Hatanaka2002,Tanjad2011}, distortion
mitigation techniques such as predistorters for power amplifier linearization~\cite{Cavers1999,Ba2008},
the approximation of nonlinear vector fields obtained from state equations~\cite{Belkhouche2005},
the obtainment of approximate solutions in simulations with complex
nonlinear systems~\cite{Storace2004}, or the search for canonical
piecewise linear representations in one and multiple dimensions~\cite{Julian1999}.

In the last decades, the main efforts in piecewise linearization have
been devoted both to find approximations of multidimensional functions
from a mathematical standpoint and to define circuit architectures
implementing them (see, for example, \cite{Brox2013} and references
therein). In the one-dimensional setting, a simple and common linearization
strategy consists in building a linear interpolant between samples
of the nonlinear function over a uniform partition of its domain.
Such a polygonal (i.e., continuous piecewise linear) interpolant may
be further optimized by choosing a better partition of the domain
according to the minimization of some error measure. This is a sensible
strategy in problems where there is a constraint on the budget of
samples allowed in the partition. 

Hence, in spite of the multiple benefits derived from modeling with
piecewise linear representations, a proper selection of the interval
partitions and/or predefining the number of partitions is paramount
for a satisfactory performance. Some researchers~\cite{Tanjad2011}
use cross-validation based approaches to select such a number of pieces
within a partition. In other applications, the budget of pieces may
be constrained by an internal design requirement (speed, memory or
target error) of the approximation algorithm or by some external condition.

Simplified models may be built using descent methods~\cite{Usow1967},
dynamic programming~\cite{Bellman1969} or heuristics such as genetic~\cite{Hatanaka2002}
and/or clustering~\cite{Ghosh2011} algorithms to optimize some target
approximation error. In some cases, however, the resulting piecewise
representation may fail to preserve desirable properties of the original
nonlinear system such as continuity~\cite{Hatanaka2002}.

We consider the simplified model representation given by the least-first-power
or best $\Lone$ approximation of a continuous nonlinear function
$f$ by some polygonal function. The generic topic of least-first-power
approximation has been previously considered in several references,
e.g., \cite{Rice1964a,Rice1964b,pinkus1989l1}, over a span of many
years and it is a recurrent topic and source of insightful results. 

We develop a fast and practical method to compute a suboptimal partition
of the interval where the polygonal interpolant and the best $\Lone$
polygonal approximation to a nonlinear function are to be computed.
This technique allows to further optimize the $\Lone$ polygonal approximation
to a function among all possible partitions having the same number
of segments, or conversely, allow to achieve a target approximation
error while minimizing the budget of segments used to represent the
nonlinear function. The resulting polygonal approximation is useful
in applications where the evaluation of continuous mathematical functions
constitutes a significant computational burden, such as computer vision~\cite{gallego2008segmentation,guillaumin2012face}
or signal processing~\cite{Xie2012,Sehili2012,Gallego2013TSP}.

Our work may be generalized to the linearization of multidimensional
functions~\cite{Julian1999} and the incorporation of constraints,
thus opening new perspectives also in the context of designing circuit
architectures for such piecewise linear approximations, as in~\cite{Brox2013}.
However, these interesting generalizations will be the topic of future
work.

The paper is organized as follows: two polygonal approximations of
real-valued univariate functions (interpolant and best $\Lone$ approximation)
are presented in Section~\ref{sec:Piecewise-Linearization}. The
mathematical foundation and algorithmic procedure to compute a suboptimal
partition for the polygonal approximations are developed in Section~\ref{sec:Optimizing-the-partition}.
The implementation of the numerical evaluation of polygonal approximations
is discussed in Section~\ref{sec:Complexity}. Experimental results
of the developed technique on nonlinear functions (Gaussian, chirp)
are given in Section~\ref{sec:Experiments}, both in terms of quality
and computational times. Finally, some conclusions are drawn in Section~\ref{sec:Conclusions}.

\section{Piecewise Linearization\label{sec:Piecewise-Linearization}}

\begin{figure}
\begin{centering}
\includegraphics[width=1\columnwidth]{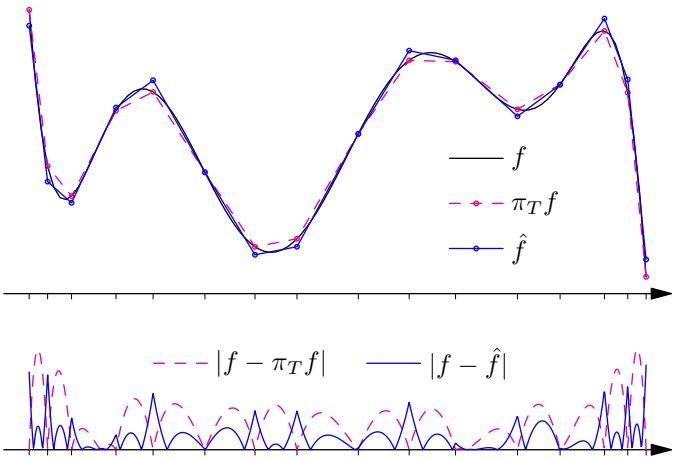}
\par\end{centering}

\caption{\label{fig:pipvsproj}Top: seventh degree polynomial $\f(x)=(x+4)(x+3)(x+2.5)x(x-1.5)(x-2)(x-3)$
and two polygonal approximations: the linear interpolant $\pip\f$
and the best $\Lone$ approximation $\fbest$. Bottom: corresponding
absolute approximation errors (magnified by a $5\times$ factor).}

\end{figure}
\label{sec:DefPiecewiseFunction}In general, a piecewise function
over an interval $I=[a,b]$ is specified by two elements: a set of
control or nodal points $\{x_{i}\}_{i=0}^{N}$, also called knots,
that determine a partition $\part=\{I_{i}\}_{i=1}^{N}$ of $I$ into
a set of $N$ (disjoint) subintervals $I_{i}=[x_{i-1},x_{i}]\mid a=x_{0}<x_{1}<\ldots<x_{N}=b$,
and a collection of $N$ functions $f_{i}(x)$ (so called ``pieces''),
one for each subinterval $I_{i}$. In particular, a \emph{polygonal}
or continuous piecewise linear (CPWL) function satisfies additional
constraints: all ``pieces'' $f_{i}(x)$ are (continuous) linear
segments and there are no jumps across pieces, i.e., continuity is
also enforced at subinterval boundaries, $f_{i}(x_{i})=f_{i+1}(x_{i})\;\forall i=\left\{ 1,\ldots,N-1\right\} $.
Fig.~\ref{fig:pipvsproj} shows, for a given partition $T,$ the
two polygonal functions that we use throughout the paper to approximate
a real-valued function $f$: the interpolant $\pip f$ and best $\Lone$
approximation~$\fbest$. Polygonal functions of a given partition
$\part$ generate a vector space $\vecspace$ since the addition of
such functions and/or multiplication by a scalar yields another polygonal
function defined over the same partition. 

\begin{figure}
\begin{centering}
\includegraphics[width=0.9\columnwidth]{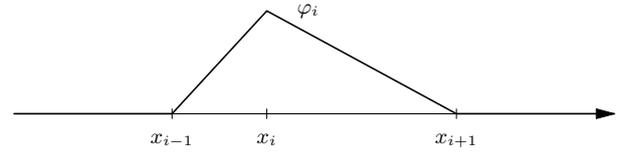}
\par\end{centering}

\caption{\label{fig:Hat-functions}Nodal basis function $\varphi_{i}$ of $\vecspace$
centered at the $i$-th control point $x_{i}$. Function $\varphi_{i}$
has the shape of a \emph{hat}; in particular, it takes value 1 at
$x_{i}$ and zero at all other control points $x_{j}$, $j\neq i$.}
\end{figure}

A useful basis for vector space $\vecspace$ is formed by the set
of nodal basis or hat functions $\left\{ \varphi_{i}\right\} _{i=0}^{N}$,
where $\varphi_{i}$, displayed in Fig.~\ref{fig:Hat-functions},
is the piecewise linear function in $\vecspace$ whose value is~1
at $x_{i}$ and zero at all other control points $x_{j}$, $j\neq i$,
i.e., 
\[
\varphi_{i}(x)=\begin{cases}
\frac{x-x_{i-1}}{x_{i}-x_{i-1}} & \text{if }x\in[x_{i-1},x_{i}],\\
\frac{x_{i+1}-x}{x_{i+1}-x_{i}} & \text{if }x\in[x_{i},x_{i+1}],\\
0 & \text{otherwise}.
\end{cases}
\]
Functions $\varphi_{0}$ and $\varphi_{N}$ associated to boundary
points $x_{0}$ and $x_{N}$ are only half hats. These basis functions
are convenient since they can represent any function $v\in\vecspace$
in terms of the values of $v$ at the control points, $v_{i}=v(x_{i})$,
in the form
\begin{equation}
v(x)=\sum_{i=0}^{N}v_{i}\varphi_{i}(x).\label{eq:ExpansionInVecSpace}
\end{equation}
From an approximation point of view this basis is frequently used
in the Finite Element Method since the hat function (\emph{simplex}
in arbitrary dimensions) is flexible, economic and in some way a natural
geometrical element into which to decompose an arbitrary geometric
object.

The polygonal interpolant $\pip f\in\vecspace$ of a continuous function
$f$ (possibly not in $\vecspace$) over the interval $I$ linearly
interpolates the samples of $f$ at the control points, thus using
$v_{i}=f(x_{i})$ in~(\ref{eq:ExpansionInVecSpace}), 
\[
\pip f(x)=\sum_{i=0}^{N}f(x_{i})\varphi_{i}(x).
\]
This polygonal approximation is trivial to construct and might be
good enough in some applications (e.g., power amplifier predistorters~\cite{Ba2008},
the trapezoidal rule for integration), but for us it is useful to
analyze other possible polygonal approximations, such as the best
one in the $\Lone$ sense, as we discuss next.

\label{sub:Best-LOne-polygonal-computation}

Now, consider the problem of approximating a continuous function $f$
by some polygonal function in~$\vecspace$ using the $\Lone$ norm
to measure distances. We address natural questions such as the existence
and uniqueness of such a best approximation, methods to determine
it and the derivation of estimates for the minimal distance.

Let us answer the question about the existence of a best approximation,
i.e., the existence of $\fbest\in\vecspace$ whose distance from $f$
is least. Recall that the space of continuous functions in a given
closed interval $I=[a,b]$, together with the $\Lone$ norm
\begin{equation}
\|u\|_{\Lone(I)}\coloneqq\int_{I}|u(x)|\,\dif x\label{eq:DefLOneNorm}
\end{equation}
is a normed linear vector space (NLVS) $(C(I),\|\cdot\|_{\Lone(I)})$.
Since $\vecspace\subset C(I)$ is a finite dimensional linear subspace
(with basis given by the nodal functions $\{\varphi_{i}\}$) of the
normed space $(C(I),\|\cdot\|_{\Lone(I)})$, then for every $\f\in C(I)$
there exists a best approximation to $\f$ in $\vecspace$ \cite[Cor. 15.10]{Plato2003}
\cite[Thm. I.1]{Rivlin1969}. 

The uniqueness of the best approximation is guaranteed for strictly
convex subspaces of NLVSs \cite[Thm. 15.19]{Plato2003} \cite[Thm. I.3]{Rivlin1969},
i.e., those whose unit balls are strictly convex sets. Linear vector
spaces with the 1 or $\infty$ norms are not strictly convex, therefore
(a priori) the solution might be unique, but it is not guaranteed.
In these cases, the uniqueness question requires special consideration.
Further insights about this topic are given in \cite[ch. 4]{Rice1964BookVol1}\cite[ch. 3]{Rivlin1969},
which are general references for $\Lone$ approximation (using polynomials
or other functions) and in \cite{DeVore1998}, which is a comprehensive
and advanced reference about nonlinear approximation. 

Next, we show how to compute such a best $\Lone$ approximation, and
later we will carry out an error analysis. As is well known~\cite[p. 130]{MoonStirling2000},
generically, the analytic approach to optimization problems using
the $\Lone$ norm involves derivatives of the absolute value, which
makes the search for an analytical solution significantly more difficult
than other problems (e.g., those using the $\Ltwo$ norm). 

As already seen, a function $v\in\vecspace$ can be written as~(\ref{eq:ExpansionInVecSpace}).
Let us explicitly note the dependence of $v$ with respect to the
coefficients $\bv=(v_{0},\ldots,v_{N})^{\top}$ by $v(x)\equiv v(x;\bv)$.
The least-first-power or best $\Lone$ approximation to $f\in C(I)$
is a function $\fbest\in\vecspace$ that minimizes $\|f-\fbest\|_{\Lone(I)}$.
Since $\fbest\in\vecspace$, it admits the expansion in terms of the
basis functions of $\vecspace$, i.e., letting $\by=(y_{0},\ldots,y_{N})^{\top}$
we may write
\begin{equation}
\fbest(x)\equiv\fbest(x;\by)=\sum_{i=0}^{N}y_{i}\varphi_{i}(x).\label{eq:DefBestLOneApproxCoeffsY}
\end{equation}
By definition, the coefficients $\by$ minimize the cost function
\begin{equation}
\cost(\mathbf{v})\coloneqq\|f(x)-\fbest(x;\bv)\|_{\Lone(I)}.\label{eq:CostFuncL1}
\end{equation}
Hence, they solve the necessary optimality conditions given by the
non-linear system of $N+1$ equations
\begin{equation}
\bg(\bv)\coloneqq\frac{\partial}{\partial\bv}\cost(\bv)=\mathbf{0},\label{eq:OptConditionsL1ComputeY}
\end{equation}
where $\bg=(g_{0},\ldots,g_{N})^{\top}$ is the gradient of~(\ref{eq:CostFuncL1}),
with entries $g_{j}(\bv)=\frac{\partial}{\partial v_{j}}\cost(\bv)$. 

Due to the partition of the interval $I$ into disjoint subintervals
$I_{i}=[x_{i-1},x_{i}]$, we may write~(\ref{eq:CostFuncL1}) as
\begin{align*}
\cost(\bv) & =\sum_{i=1}^{N}\|f(x)-\fbest(x;\bv)\|_{\Lone(I_{i})}\\
 & =\sum_{i=1}^{N}\int_{I_{i}}\left|\f(x)-\left(v_{i-1}\varphi_{i-1}(x)+v_{i}\varphi_{i}(x)\right)\right|\dif x,
\end{align*}
therefore, if $\signum(x)=x/|x|$, each gradient component is,
\begin{align*}
g_{j}(\bv)= & -\int_{I_{j}}\signum\Bigl(\f(x)-\sum_{n=j-1}^{j}v_{n}\varphi_{n}(x)\Bigr)\varphi_{j}(x)\,\dif x\\
 & -\int_{I_{j+1}}\signum\Bigl(\f(x)-\sum_{n=j}^{j+1}v_{n}\varphi_{n}(x)\Bigr)\varphi_{j}(x)\,\dif x.
\end{align*}
Observe that $g_{j}$ solely depends on $\{v_{j-1},v_{j},v_{j+1}\}$
(except at extreme cases $j=\{0,N\}$) due to the locality and adjacency
of the basis functions $\{\varphi_{i}\}$. In the case of the $\Ltwo$
norm, the optimality conditions are linear and the previous observation
leads to a tridiagonal (linear) system of equations. 

A closed form solution of~(\ref{eq:OptConditionsL1ComputeY}) may
not be available, and so, to solve the system we use standard numerical
iterative algorithms of the form $\bv^{k+1}=\bv^{k}+\step^{k}$ to
find an approximate local solution $\by=\lim_{k\to\infty}\bv^{k}$.
Specifically, we use step $\step^{k}$ in the Newton-Raphson iteration
given by the solution of the linear system $\mH(\bv^{k})\,\step^{k}=-\bg(\bv^{k})$,
where $\mH=\frac{\partial\bg}{\partial\bv}$. Near the solution~$\by$,
this iteration has quadratic convergence rate. This is also the step
given by Newton's optimization method when approximating~(\ref{eq:CostFuncL1})
by its quadratic Taylor model. Due to the locality of the basis functions
$\{\varphi_{i}\}$, cost function~(\ref{eq:CostFuncL1}) has the
advantage that its Hessian ($\mH$) is a tridiagonal matrix, so $\step^{k}$
is faster to compute than in the case of a full Hessian matrix. 

The search for the optimal coefficients may be initialized by setting
$\bv^{0}$ equal to the values of the function at the nodal points,
i.e., $\bv^{0}=(v_{0}^{0},\ldots,v_{N}^{0})^{\top}$ with $v_{i}^{0}=f(x_{i})$.
A more sensible initialization $\bv^{0}$ to improve convergence toward
the optimal coefficients is given by the ordinates $\bv^{0}=\mathbf{c}=(c_{0},\ldots,c_{N})^{\top}$
of the best $\Ltwo$ approximation (i.e., orthogonal projection of
$f$ onto $\vecspace$) $\ortproj\f(x)=\sum_{i=0}^{N}c_{i}\varphi_{i}(x)$,
which are easily obtained by solving a linear system of equations
using the Thomas algorithm, $M\mathbf{c}=\mathbf{b}$ with tridiagonal
Gramian matrix $M=(m_{ij}),$ $m_{ij}=\inner{\varphi_{i}}{\varphi_{j}}$,
$\mathbf{b}=(b_{0},\ldots,b_{N})^{\top}$, $b_{i}=\inner f{\varphi_{i}}$
and inner product $\inner uv\coloneqq\int_{I}u(x)v(x)\dif x$.

From a numerical point of view, it is also a reasonable choice to
replace $\signum(x)$ by some smooth approximation, for example, $\signum(x)\approx\tanh(kx)$,
with parameter $k\gg1$ controlling the width of the transition around
$x=0$. 

In summary, the coefficients $\by$ that specify the best $\Lone$
approximation~(\ref{eq:DefBestLOneApproxCoeffsY}) on a given partition
$\part$ are computed numerically via iterative local optimization
techniques starting from an initial guess $\bv^{0}$.

\section{Optimizing the partition\label{sec:Optimizing-the-partition}}

Given a vector space $\vecspace$, we are endowed with a procedure
to compute the least-first-power approximation of a function $f$
and the corresponding error, $\|f-\fbest\|_{\Lone(I)}$. However,
the approximation error depends on the choice of $\vecspace$, which
is specified by the partition $\part$. Hence, the next problem that
naturally arises is the optimization of the partition $\part$ for
a given budget of control points, i.e., the search for the best vector
space $\vecspace$ to approximate $f$ for a given partition size.
This is a challenging non-linear optimization problem, even in the
simpler case (less degrees of freedom) of substituting $\fbest$ by
the polygonal interpolant~$\pip f$. Fortunately, a good approximation
of the optimal partition $\partOpt$ can be easily found using an
asymptotic analysis.

Next, we carry out a detailed error analysis for the polygonal interpolant
$\pip f$ and the polygonal least-first-power approximation $\fbest$.
This will help us derive an approximation to the optimal partition
that is valid for both $\pip f$ and $\fbest$, because, as it will
be shown, their approximation errors are roughly proportional if a
sufficiently large budget of control points, i.e., large number of
subintervals, is available.

\begin{figure}
\begin{centering}
\includegraphics[width=1\columnwidth]{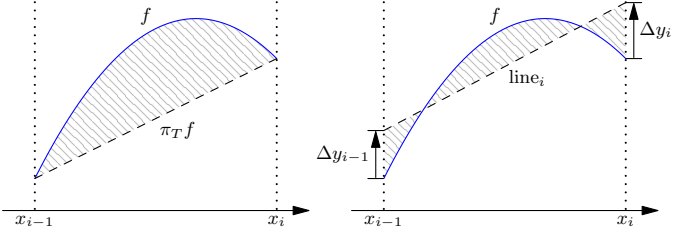}
\par\end{centering}

\caption{\label{fig:errorzone}Function $f$ and two linear approximations
in interval $I_{i}=[x_{i-1},x_{i}]$. Left: interpolant $\pip f$
defined in~(\ref{eq:LinearInterpolantInInterval}). Right: arbitrary
linear segment $\mathrm{line}_{i}$ defined in~(\ref{eq:DefLineSegment}),
where $\Delta y_{j}$ is a signed vertical displacement with respect
to $f(x_{j})$.}

\end{figure}

\subsection{Error in a single interval: linear interpolant\label{sub:Error-Single-Interval}}

First, let us analyze the error generated when approximating a function
$f$, twice continuously differentiable, by its polygonal interpolant
$\pip f$ in a single interval $I_{i}=[x_{i-1},x_{i}]$, of length
$h_{i}=x_{i}-x_{i-1}$ . To this end, recall the following theorem
on interpolation errors~\cite[sec. 4.2]{Cheney2012}: Let $f$ be
a function in $C^{n+1}(\Omega)$, with $\Omega\subset\R$ and closed,
and let $p$ be a polynomial of degree $n$ or less that interpolates
$f$ at $n+1$ distinct points $x_{0},\ldots,x_{n}\in\Omega$. Then,
for each $x\in\Omega$ there exists a point $\xi_{x}\in\Omega$ for
which 
\begin{equation}
\f(x)-p(x)=\frac{1}{(n+1)!}\f^{(n+1)}(\xi_{x})\prod_{i=0}^{n}(x-x_{i}).\label{eq:ThmErrorPolInterp}
\end{equation}
In the subinterval $I_{i}$, letting $\delta_{i}(x)=(x-x_{i-1})/h_{i}$,
the polygonal interpolant $\pip f$ is written as
\begin{equation}
\pip\f(x)=\f(x_{i-1})\bigl(1-\delta_{i}(x)\bigr)+\f(x_{i})\delta_{i}(x).\label{eq:LinearInterpolantInInterval}
\end{equation}
Since $\pip\f$ interpolates the function $\f$ at the endpoints of
$I_{i}$, we can apply theorem~(\ref{eq:ThmErrorPolInterp}) (with
$n=1$); hence, the approximation error only depends on $\fsec$ and
$x$, but not on $\f$ or $\f^{\prime}$:
\begin{equation}
\f(x)-\pip\f(x)=-\frac{1}{2}\fsec(\xi_{x})(x-x_{i-1})(x_{i}-x).\label{eq:interpolant-difference}
\end{equation}
Let us compute the $\Lone$ error over the subinterval $I_{i}$ by
integrating the magnitude of~(\ref{eq:interpolant-difference}),
according to~(\ref{eq:DefLOneNorm}):
\begin{align*}
\|\f-\pip\f\|_{\Lone(I_{i})} & =\int_{I_{i}}\left|-\frac{1}{2}\fsec(\xi_{x})(x-x_{i-1})(x_{i}-x)\right|\,\dif x\\
 & =\frac{1}{2}\int_{I_{i}}\left|\fsec(\xi_{x})\right|\left|(x-x_{i-1})(x_{i}-x)\right|\dif x.
\end{align*}
Next, recall the first mean value theorem for integration, which states
that if $u:[A,B]\rightarrow\R$ is a continuous function and $v$
is an integrable function that does not change sign on $(A,B)$, then
there exists a number $\xi\in(A,B)$ such that
\begin{equation}
\int_{A}^{B}u(x)v(x)\,\dif x=u(\xi)\int_{A}^{B}v(x)\,\dif x.\label{eq:FirstMVTIntegration}
\end{equation}
Applying~(\ref{eq:FirstMVTIntegration}) to the $\Lone$ error and
noting that $(x-x_{i-1})(x_{i}-x)\geq0$ $\forall x\in I_{i}$ gives
\begin{align}
\|\f-\pip\f\|_{\Lone(I_{i})} & \stackrel{\eqref{eq:FirstMVTIntegration}}{=}\frac{1}{2}\left|\fsec(\eta)\right|\int_{I_{i}}(x-x_{i-1})(x_{i}-x)\,\dif x\nonumber \\
 & =\frac{h_{i}^{3}}{12}\left|\fsec(\eta)\right|,\label{eq:interpolant-error-in-an-interval}
\end{align}
for $\eta\in(x_{i-1},x_{i})$. Finally, if $\fimax\coloneqq\max_{\eta\in I_{i}}|\fsec(\eta)|$,
a direct derivation of the $\Lone$ error bound yields
\begin{equation}
\|\f-\pip\f\|_{\Lone(I_{i})}\leq\frac{1}{12}\fimax h_{i}^{3}.\label{eq:interpolation-bound-for-L1error-in-an-interval}
\end{equation}

Formula~(\ref{eq:interpolation-bound-for-L1error-in-an-interval})
states that the deviation of $f$ from being linear between endpoints
of $I_{i}$ is bounded by the maximum concavity/convexity of the function
in $I_{i}$ (e.g., $\fimax$ limits the amount of bending) and the
cubic power of the interval size $h_{i}$, also known as the local
density of control points.

\subsection{Error in a single interval: best $\Lone$ linear approximation\label{sub:Minimum-error-Line-segment}}

To analyze the error due to the least-first-power approximation $\fbest$
and see how much it improves over that of the interpolant $\pip f$,
let us first characterize the error incurred when approximating a
function $\f(x)$ by a linear segment not necessarily passing through
the endpoints of $I_{i}$,
\begin{align}
\linei(x;\Delta y_{i-1},\Delta y_{i}) & =\bigl(\f(x_{i-1})+\Delta y_{i-1}\bigr)\bigl(1-\delta_{i}(x)\bigr)\nonumber \\
 & \quad+\bigl(\f(x_{i})+\Delta y_{i}\bigr)\delta_{i}(x),\label{eq:DefLineSegment}
\end{align}
where $\Delta y_{i-1}$ and $\Delta y_{i}$ are extra parameters with
respect to $\pip\f$ that allow the linear segment to better fit the
function $\f$ in $I_{i}$. Letting $(\pip\Delta y)(x;\Delta y_{i-1},\Delta y_{i})\coloneqq\Delta y_{i-1}\bigl(1-\delta_{i}(x)\bigr)+\Delta y_{i}\delta_{i}(x)$
by analogy to~(\ref{eq:LinearInterpolantInInterval}), the corresponding
error $\epsilon\equiv\epsilon(x;\Delta y_{i-1},\Delta y_{i})$ is
\begin{align}
\epsilon & =\f(x)-\linei(x;\Delta y_{i-1},\Delta y_{i}),\label{eq:DefEpsilonLineSegment}\\
 & =\f(x)-\pip\f(x)-(\pip\Delta y)(x;\Delta y_{i-1},\Delta y_{i}),\nonumber \\
 & \stackrel{\eqref{eq:interpolant-difference}}{=}-\frac{1}{2}\fsec(\xi_{x})(x-x_{i-1})(x_{i}-x)\nonumber \\
 & \quad\;-(\pip\Delta y)(x;\Delta y_{i-1},\Delta y_{i}).\label{eq:EpsilonSimplified}
\end{align}

\subsubsection{Characterization of the optimal line segment}

To find the line segment that minimizes the $\Lone$ distance
\[
\|\epsilon\|_{\Lone(I_{i})}=\|\f-\linei\|_{\Lone(I_{i})}=\int_{I_{i}}\left|\epsilon(x;\Delta y_{i-1},\Delta y_{i})\right|\dif x,
\]
i.e., to specify the values of the optimal $\Delta y_{i-1},\Delta y_{i}$
in~(\ref{eq:DefLineSegment}), we solve the necessary optimality
conditions given by the non-linear system of equations
\begin{equation}
\begin{aligned}0=\frac{\partial\|\epsilon\|_{\Lone(I_{i})}}{\partial\Delta y_{i-1}} & =\int_{I_{i}}\signum\bigl(\epsilon(x;\Delta y_{i-1},\Delta y_{i})\bigr)\bigl(1-\delta_{i}(x)\bigr)\,\dif x,\\
0=\frac{\partial\|\epsilon\|_{\Lone(I_{i})}}{\partial\Delta y_{i}} & =\int_{I_{i}}\signum\bigl(\epsilon(x;\Delta y_{i-1},\Delta y_{i})\bigr)\delta_{i}(x)\,\dif x,
\end{aligned}
\label{eq:OptConditionsBestSegmentL1}
\end{equation}
where we used that, for a function $g(x)$,
\[
\frac{\partial}{\partial x}|g(x)|=\frac{\partial}{\partial x}\sqrt{g^{2}(x)}=\signum\bigl(g(x)\bigr)\frac{\partial}{\partial x}g(x).
\]

Adding both optimality equations in~(\ref{eq:OptConditionsBestSegmentL1})
gives
\[
\int_{I_{i}}\signum\bigl(\epsilon(x;\Delta y_{i-1},\Delta y_{i})\bigr)\dif x=0,
\]
which implies that $\epsilon$ must be positive in half of the interval
$I_{i}$ and negative in the other half.

In fact, \cite{KripkeRivlin1965}\cite[Cor. 3.1.1]{Rivlin1969} state
that if $\epsilon$ has a finite number of zeros (at which $\epsilon$
changes sign) in $I_{i}$, then $\linei$ is a best $\Lone$ approximation
to $\f$ if and only if~(\ref{eq:OptConditionsBestSegmentL1}) is
satisfied. To answer the uniqueness question, \cite{Jackson1921Note}\cite{Jackson1930}\cite[Thm 3.2]{Rivlin1969}
state that a continuous function on $I_{i}$ has a unique best $\Lone$
approximation out of the set of polynomials of degree $\leq n$. Hence,
the solution of~(\ref{eq:OptConditionsBestSegmentL1}) provides the
best $\Lone$ linear approximation.

Let us discuss the solution of~(\ref{eq:OptConditionsBestSegmentL1}).
If $\epsilon$ changes sign only at one abscissa $\bar{x}\in I_{i}$,
e.g., $\epsilon(\bar{x})=0$, $\epsilon(\{x<\bar{x}\})<0$ and $\epsilon(\{x>\bar{x}\})>0$,
the non-linear system of equations~(\ref{eq:OptConditionsBestSegmentL1})
cannot be satisfied since the first equation gives $\bar{x}=x_{i-1}+h_{i}(1-1/\sqrt{2})$
while the second equation gives $\bar{x}=x_{i-1}+h_{i}/\sqrt{2}$.
However, in the next simplest case where $\epsilon$ changes sign
at two abscissas $\bar{x}_{1},\bar{x}_{2}\in I_{i}$, the non-linear
system~(\ref{eq:OptConditionsBestSegmentL1}) does admit a solution.
This is also intuitive to justify since it corresponds to the simplified
case $\fsec=C$ constant in~$I_{i}$, where the sign change occurs
if $\epsilon=0$, i.e., according to~(\ref{eq:EpsilonSimplified}),
$\frac{1}{2}C\,(x-x_{i-1})(x_{i}-x)+\Delta y_{i-1}\bigl(1-\delta_{i}(x)\bigr)+\Delta y_{i}\delta_{i}(x)=0,$
which is a quadratic equation in $x$. It is also intuitive by looking
at a plot of a candidate small error, as in Fig.~\ref{fig:errorzone},
right.

Next, we further analyze the aforementioned case of $\epsilon$ changing
sign at $\bar{x}_{1},\bar{x}_{2}\in I_{i}$, with $\bar{x}_{2}>\bar{x}_{1}$.
Assume that $\signum(\epsilon)=-1$ for $\bar{x}_{1}<x<\bar{x}_{2}$
and $\signum(\epsilon)=+1$ in the other half of $I_{i}$. If we apply
the change of variables $t=\delta_{i}(x)=(x-x_{i-1})/h_{i}$, and
let $t_{j}=\delta_{i}(\bar{x}_{j})$ for $j=1,2$, then~(\ref{eq:OptConditionsBestSegmentL1})
becomes 
\begin{align*}
t_{2}^{2}-t_{1}^{2}-2(t_{2}-t_{1})+\frac{1}{2} & =0,\\
t_{1}^{2}-t_{2}^{2}+\frac{1}{2} & =0.
\end{align*}
 Adding both equations gives, as we already mentioned, $t_{2}-t_{1}=\frac{1}{2}$,
i.e., $\bar{x}_{2}-\bar{x}_{1}=\frac{1}{2}h_{i}$, stating that $\epsilon<0$
in half $ $of the interval. This equation can be used to simplify
the second equation, $(t_{2}+t_{1})(t_{2}-t_{1})=\frac{1}{2}$, yielding
$t_{2}+t_{1}=1$. Therefore~(\ref{eq:OptConditionsBestSegmentL1})
is equivalent to the linear system $\{t_{2}-t_{1}=\frac{1}{2},\; t_{2}+t_{1}=1\},$
whose solution is $t_{1}=\frac{1}{4}$, $t_{2}=\frac{3}{4}$, i.e.
$\bar{x}_{1}=x_{i-1}+\frac{1}{4}h_{i}$, $\bar{x}_{2}=x_{i-1}+\frac{3}{4}h_{i}$. 

This agrees with the particularization of a more general result \cite[Cor.3.4.1]{Rivlin1969}:
if $\f$ is adjoined to the set of (linear, $n=1$) polynomials in
$I_{i}$, $P_{n}(I_{i})$, meaning that $f\in C(I_{i})\setminus P_{n}(I_{i})$
and $f-p$ has at most $n+1$ distinct zeros in $I_{i}$ for every
$p\in P_{n}(I_{i})$, its best $\Lone$ approximation out of $P_{n}(I_{i})$
is the unique $\ell^{\ast}\in P_{n}(I_{i})$ which satisfies
\[
\ell^{\ast}\left(\bar{x}_{j}\right)=\f\left(\bar{x}_{j}\right)
\]
 for $\bar{x}_{j}=x_{i-1}+\bigl(1+\cos(j\pi/(n+2))\bigr)h_{i}/2\in I_{i}$,
$j=1,\ldots,n+1$. The cosine term comes from the zeros of the Chebyshev
polynomial of the second kind. 

In other words, the best approximation is constructed by interpolating
$f$ at the \emph{canonical points} $\bar{x}_{j}$ (the points of
sign change of $\signum(\text{\ensuremath{\epsilon}})$ in~(\ref{eq:OptConditionsBestSegmentL1})),
as expressed by~\cite{Rice1964b}\cite{Rice1964BookVol1} in a nonlinear
context. Hence, the values of $\Delta y_{i-1},\Delta y_{i}$ that
satisfy~(\ref{eq:OptConditionsBestSegmentL1}) are chosen so that
zero crossings of $\epsilon(x;\Delta y_{i-1},\Delta y_{i})$ occur
at canonical points $\frac{1}{4}$ and $\frac{3}{4}$ length of the
interval $I_{i}$, yielding the linear system of equations
\[
\left.\begin{array}{c}
\epsilon(\bar{x}_{1};\Delta y_{i-1},\Delta y_{i})=0\\
\epsilon(\bar{x}_{2};\Delta y_{i-1},\Delta y_{i})=0
\end{array}\right\} 
\]
 whose solution is, after substituting in~(\ref{eq:EpsilonSimplified}),
\begin{equation}
\begin{aligned}\Delta y_{i-1} & =\frac{3h_{i}^{2}}{64}\left(\fsec(\xi_{\bar{x}_{2}})-3\fsec(\xi_{\bar{x}_{1}})\right),\\
\Delta y_{i} & =\frac{3h_{i}^{2}}{64}\left(-3\fsec(\xi_{\bar{x}_{2}})+\fsec(\xi_{\bar{x}_{1}})\right).
\end{aligned}
\label{eq:SolutionDeltaY}
\end{equation}

The previous solution implies that the sum of the displacements has
opposite sign to the convexity/concavity of the function $\f$:
\begin{equation}
\Delta y_{i-1}+\Delta y_{i}=-\frac{3h_{i}^{2}}{16}\fsec(\eta),\label{eq:SumDeltays}
\end{equation}
where $\fsec(\xi_{\bar{x}_{1}})+\fsec(\xi_{\bar{x}_{2}})$ from~(\ref{eq:SolutionDeltaY})
lies between the least and greatest values of $2\fsec$ on $I_{i}$,
and by the intermediate value theorem it is $2\fsec(\eta)$ for some
$\eta\in(x_{i-1},x_{i})$. This agrees with the intuition/graphical
interpretation (see Fig.~\ref{fig:errorzone}, right).

\subsubsection{Minimum error of the optimal line segment}

Now that the optimal $\Delta y_{i-1},\Delta y_{i}$ have been specified,
we may compute the minimum error. Let $s=\signum\bigl(\epsilon(x;\Delta y_{i-1},\Delta y_{i})\bigr)=\pm1$
for $\bar{x}_{1}<x<\bar{x}_{2}$, then, since $|a|=\signum(a)a$,
we may expand
\begin{equation}
\min\|\epsilon\|_{\Lone(I_{i})}=\left(-\int_{x_{i-1}}^{\bar{x}_{1}}\epsilon\,\dif x+\int_{\bar{x}_{1}}^{\bar{x}_{2}}\epsilon\,\dif x-\int_{\bar{x}_{2}}^{x_{i}}\epsilon\,\dif x\right)s.\label{eq:miniCostL1-tmp}
\end{equation}

Next, since $(x-x_{i-1})(x_{i}-x)\geq0$ for all $x\in[p,q]\subset I_{i}$,
use the first mean value theorem for integration~(\ref{eq:FirstMVTIntegration})
to simplify
\begin{eqnarray*}
\int_{p}^{q}\epsilon\,\dif x & \!\!\!\!\stackrel{\eqref{eq:EpsilonSimplified}\eqref{eq:FirstMVTIntegration}}{=}\!\!\!\! & -\frac{1}{2}\fsec(\eta_{pq})\int_{p}^{q}(x-x_{i-1})(x_{i}-x)\,\dif x\\
 &  & -\Delta y_{i-1}\int_{p}^{q}\bigl(1-\delta_{i}(x)\bigr)\,\dif x-\Delta y_{i}\int_{p}^{q}\delta_{i}(x)\,\dif x\\
 & = & -\frac{1}{2}\fsec(\eta_{pq})h_{i}^{3}\left[\frac{\delta_{i}^{2}(x)}{2}-\frac{\delta_{i}^{3}(x)}{3}\right]_{p}^{q}\\
 &  & -\Delta y_{i-1}h_{i}\left[\delta_{i}(x)-\frac{\delta_{i}^{2}(x)}{2}\right]_{p}^{q}-\Delta y_{i}\left[\frac{\delta_{i}^{2}(x)}{2}\right]_{p}^{q},
\end{eqnarray*}
for some $\eta_{pq}\in(p,q)$. In particular, using the previous formula
for each term in~(\ref{eq:miniCostL1-tmp}) gives
\begin{align*}
-\int_{x_{i-1}}^{\bar{x}_{1}}\epsilon\,\dif x & =\frac{1}{2}\fsec(\eta_{1})\frac{5h_{i}^{3}}{192}+\Delta y_{i-1}\frac{7h_{i}}{32}+\Delta y_{i}\frac{h_{i}}{32},\\
\int_{\bar{x}_{1}}^{\bar{x}_{2}}\epsilon\,\dif x & =-\frac{1}{2}\fsec(\eta_{2})\frac{22h_{i}^{3}}{192}-\Delta y_{i-1}\frac{8h_{i}}{32}-\Delta y_{i}\frac{8h_{i}}{32},\\
-\int_{\bar{x}_{2}}^{x_{i}}\epsilon\,\dif x & =\frac{1}{2}\fsec(\eta_{3})\frac{5h_{i}^{3}}{192}+\Delta y_{i-1}\frac{h_{i}}{32}+\Delta y_{i}\frac{7h_{i}}{32},
\end{align*}
for some $\{\eta_{1},\eta_{2},\eta_{3}\}\in(x_{i-1},x_{i})$. Hence,
(\ref{eq:miniCostL1-tmp}) becomes
\begin{align*}
\min\|\epsilon\|_{\Lone(I_{i})} & =s\frac{h_{i}^{3}}{384}\left(5\fsec(\eta_{1})-22\fsec(\eta_{2})+5\fsec(\eta_{3})\right).
\end{align*}

The segments in the best polygonal $\Lone$ approximation $\fbest$
may not strictly satisfy this because $\fbest$ has additional continuity
constraints across segments. The jump discontinuity at $x=x_{i}$
between adjacent independently-optimized pieces is
\begin{align*}
\left|\rightDeltaY i^{-}-\leftDeltaY i^{+}\right| & \leq\frac{3h_{i}^{2}}{64}\,\left|-3\fsec(\xi_{\bar{x}_{2},i})+\fsec(\xi_{\bar{x}_{1},i})\right|\\
 & \quad\;+\frac{3h_{i+1}^{2}}{64}\,\left|\fsec(\xi_{\bar{x}_{2},i+1})-3\fsec(\xi_{\bar{x}_{1},i+1})\right|,
\end{align*}
where $\Delta y_{i}^{-}=\left(-3\fsec(\xi_{\bar{x}_{2},i})+\fsec(\xi_{\bar{x}_{1},i})\right)3h_{i}^{2}/64$
and $\Delta y_{i}^{+}=\left(\fsec(\xi_{\bar{x}_{2},i+1})-3\fsec(\xi_{\bar{x}_{1},i+1})\right)3h_{i+1}^{2}/64$
are displacements with respect to $\f(x_{i})$ of the optimized segments~(\ref{eq:SolutionDeltaY})
at each side of $x=x_{i}$, and evaluation points $\xi_{\bar{x}_{1},j}$
and $\xi_{\bar{x}_{2},j}$ lie in $I_{j}$. In case of twice continuously
differentiable functions in a closed interval, the extreme value theorem
states that the absolute value terms in the previous equation are
bounded. Accordingly, if $h_{i}$ and $h_{i+1}$ decrease due to a
finer partition $\part$ of the interval $I$ (i.e., a larger number
of segments $N$ in $\part$), the discontinuity jumps at the control
points of the partition decrease, too. Therefore, the approximation
$\|\f-\fbest\|_{\Lone(I_{i})}\approx\min\|\f-\mathrm{line}_{i}\|_{\Lone(I_{i})}^{2}$
is valid for large $N$. Finally, if $I_{i}$ is sufficiently small
so that $\fsec$ is approximately constant within it, say $\fsec_{I_{i}}$,
then
\begin{equation}
\min\|\f-\mathrm{line}_{i}\|_{\Lone(I_{i})}\approx h_{i}^{3}s\fsec_{I_{i}}\frac{(-12)}{384}=\frac{h_{i}^{3}}{32}\left|\fsec_{I_{i}}\right|.\label{eq:approximate-best-L1-segment}
\end{equation}
In the last step we substituted $s=-\signum(\fsec_{I_{i}})$, which
can be proven by evaluation at the midpoint of interval $I_{i}$:
\begin{eqnarray*}
s & = & \signum\left(\epsilon\left(\frac{x_{i-1}+x_{i}}{2};\Delta y_{i-1},\Delta y_{i}\right)\right)\\
 & \stackrel{\eqref{eq:EpsilonSimplified}}{=} & \signum\left(-\fsec_{I_{i}}\frac{h_{i}^{2}}{4}-(\Delta y_{i-1}+\Delta y_{i})\right)\\
 & \stackrel{\eqref{eq:SumDeltays}}{=} & \signum\left(-\frac{4h_{i}^{2}}{16}\fsec_{I_{i}}+\frac{3h_{i}^{2}}{16}\fsec_{I_{i}}\right)\\
 & = & -\signum\left(\fsec_{I_{i}}\right).
\end{eqnarray*}

In the same asymptotic situation, the error of the linear interpolant~(\ref{eq:interpolant-error-in-an-interval})
becomes
\begin{equation}
\|\f-\pip\f\|_{\Lone(I_{i})}\approx\frac{h_{i}^{3}}{12}\left|\fsec_{I_{i}}\right|,\label{eq:approximate-L1-interpolated-segment}
\end{equation}
which is larger than the best $\Lone$ approximation error~(\ref{eq:approximate-best-L1-segment})
by a factor of $8/3\approx2.67$.

\subsection{\label{sub:Approximation-to-the-optimal-partition}Approximation
to the optimal partition}

Once analyzed the errors of both interpolant and least-first-order
approximation on a subinterval $I_{i}$, let us use such results to
propose a suboptimal partition $\partOpt$ of the interval $I$ in
the asymptotic case of a large number $N$ of subintervals.

A suboptimal partition for a given budget of control points $(N+1)$
is one in which every subinterval has approximately equal contribution
to the total approximation error~\cite{deBoor2001,Cox2001}. Since
such an error depends on the function $f$ being approximated, it
is clear that such a dependence will be transferred to the suboptimal
partition, i.e., the suboptimal partition is tailored to $f$. Specifically,
because the error is proportional to the local amount of convexity/concavity
of the function, a suboptimal partition places more controls points
in regions of $\f$ with larger convexity than in other regions so
that error equalization is achieved. Assuming $N$ is large enough
so that $\fsec$ is approximately constant in each subinterval and
therefore the bound~(\ref{eq:interpolation-bound-for-L1error-in-an-interval})
is tight, we have
\begin{equation}
\fimax h_{i}^{3}\approx C,\label{eq:ErrorEqualization}
\end{equation}
for some constant $C>0$, and the control points should be chosen
so that the local knot spacing~\cite{Cox2001} is $h_{i}\propto\fimax^{-1/3}$,
i.e., smaller intervals as $\fimax$ increases. Hence, the local knot
distribution or density is 
\begin{equation}
lkd(x)\propto|\fsec(x)|^{1/3},\label{eq:limit-optimal-vertex-density}
\end{equation}
meaning, as already announced, that more knots of the partition are
placed in the regions with larger magnitude of the second derivative.
\begin{figure}
\centering{}\includegraphics[width=1\columnwidth]{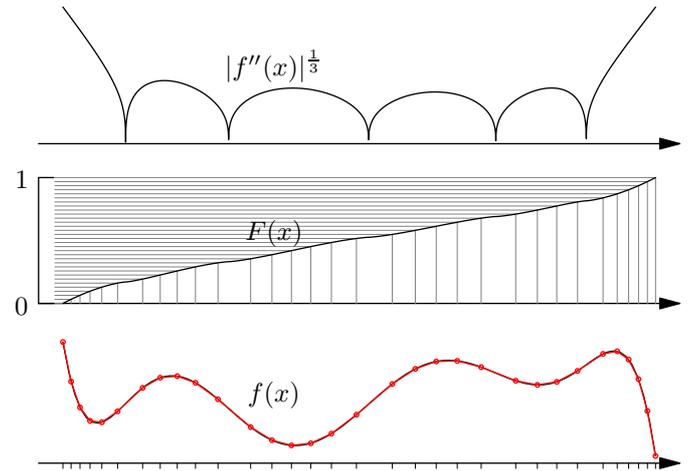}\caption{\textbf{\label{fig:GraphicalSummary}}Graphical summary of the proposed
suboptimal partition computation tailored to a given function $f$.
Top: local knot density~(\ref{eq:limit-optimal-vertex-density})
obtained from input function $f$ (Fig.~\ref{fig:pipvsproj}). Middle:
cumulative knot distribution function $F$ given by~(\ref{eq:cumulative-knot-distribution})
and control points (i.e. knots) given by the preimages of the endpoints
corresponding to a uniform partition of the range of $F$, as expressed
by~(\ref{eq:optimized-partition}). Bottom: polygonal interpolant
$\pipOpt f$ with $N=31$ ($32$ knots) overlaid on the input function
$f$. Knots are distributed according to the amount of local convexity/concavity
of $f$ displayed in top plot so that error equalization is achieved.
Hence, fewer knots are placed around the zeros of the $lkd$, which
correspond to the horizontal regions of $F$.}
\end{figure}

The error equalization criterion leads to the following suboptimal
partition $\partQuasiOpt$: $x_{0}=a$, $x_{N}=b$, and take knots
$\{x_{i}\}_{i=1}^{N-1}$ given by 
\begin{equation}
F(x_{i})=i/N,\label{eq:optimized-partition}
\end{equation}
where the monotonically increasing function $F:[a,b]\to[0,1]$ is
\begin{equation}
F(x)=\int_{a}^{x}\left|\fsec(t)\right|^{1/3}\mathrm{d}t\Bigl/\int_{a}^{b}\left|\fsec(t)\right|^{1/3}\mathrm{d}t.\label{eq:cumulative-knot-distribution}
\end{equation}
This procedure divides the range of $F(x)$ into $N$ contiguous equal
length sub-ranges, and chooses the control points $x_{i}$ given by
the preimages of the endpoints of the sub-ranges. It is graphically
illustrated in Fig.~\ref{fig:GraphicalSummary}. The suboptimal partition
is related to the theory of optimum quantization~\cite{Gersho1978},
particularly in the asymptotic or high-resolution quantization case~\cite{Lookabaugh1989},
where a ``companding'' function such as $F(x)$ enables non-uniform
subinterval spacing within a partition.

This partition allows us to estimate the error bound $\|\f-\pipQuasiOpt\f\|_{\Lone(I)}$
in the entire interval $I=[a,b]$ starting from that of the subintervals.
For any partition $\part$, the total error is the sum of the errors
over all subintervals $I_{i}$ and, by~(\ref{eq:interpolation-bound-for-L1error-in-an-interval}),
\begin{align}
\|\f-\pip\f\|_{\Lone(I)} & \leq\sum_{i=1}^{N}\frac{1}{12}\fimax h_{i}^{3},\label{eq:BoundpipL1-prev}
\end{align}
which, under the $\partQuasiOpt$ error equalization condition~(\ref{eq:ErrorEqualization}),
becomes 
\begin{equation}
\|\f-\pipQuasiOpt\f\|_{\Lone(I)}\leq\sum_{i=1}^{N}\frac{1}{12}C=\frac{1}{12}CN.\label{eq:BoundpipL1Interval}
\end{equation}
To determine $C$, sum $\fimax^{1/3}h_{i}\approx C^{1/3}$ over all
subintervals $I_{i}$ and approximate the result using the Riemann
integral: 
\begin{align}
C^{1/3}N & \approx\sum_{i=1}^{N}\fimax^{1/3}h_{i}\approx\int_{a}^{b}\left|\fsec(t)\right|^{1/3}\mathrm{d}t,\label{eq:CalculationOfConstantC-L1}
\end{align}
whose right hand side does not depend on~$N$. Substituting~(\ref{eq:CalculationOfConstantC-L1})
in~(\ref{eq:BoundpipL1Interval}) gives the approximate error bound
for the polygonal interpolant over the entire interval $I=[a,b]$:
\begin{equation}
\|\f-\pipQuasiOpt\f\|_{\Lone(I)}\lesssim\frac{1}{12N^{2}}\left(\int_{a}^{b}\left|\fsec(t)\right|^{1/3}\mathrm{d}t\right)^{3}.\label{eq:ApproxL1Boundpip}
\end{equation}

Finally, in the asymptotic case of large $N$, the approximation error
of $\fbest$ is roughly proportional to that of $\pip\f$ as shown
in~(\ref{eq:approximate-best-L1-segment}) and~(\ref{eq:approximate-L1-interpolated-segment}).
Hence, the partition specified by~(\ref{eq:optimized-partition})
is also a remarkable approximation to the optimal partition for $\fbest$
as $N$ increases. This, together with~(\ref{eq:ApproxL1Boundpip})
implies that both polygonal approximations converge to $\f$ at a
rate of at least $O(N^{-2})$.

Following a similar procedure, it is possible to estimate an error
bound on the uniform partition $\partUnif$, which can be compared
to that of the optimized one. For $\partUnif$, substitute $h_{i}=(b-a)/N$
in~(\ref{eq:BoundpipL1-prev}) and approximate the result using the
Riemann integral,
\begin{align}
\|\f-\pipUnif\f\|_{\Lone(I)} & \leq\frac{h_{i}^{2}}{12}\sum_{i=1}^{N}\fimax h_{i}\nonumber \\
 & \approx\frac{\left(b-a\right)^{2}}{12N^{2}}\|\fsec\|_{\Lone(I)}.\label{eq:ApproxL1boundUniformPip}
\end{align}
The quotient of~(\ref{eq:ApproxL1Boundpip}) and~(\ref{eq:ApproxL1boundUniformPip})
provides an estimate of the gain obtained by optimizing a partition.

\subsection{Extension to vector valued functions\label{sub:Extension-to-vectorfuncs}}

Let us point out, using a geometric approach, how the previous method
may be extended to handle vector valued nonlinear functions, i.e.,
functions with multiple values for the same $x$. Let $\vf(x)=\bigl(f_{1}(x),\ldots,f_{n}(x)\bigr)^{\top}\in\R^{n}$
consist of several functions $f_{j}(x)$ defined over the same interval
$I$, and let $\|\bv\|_{1}=\sum_{j=\text{1}}^{n}|v_{j}|$ be the usual
1-norm in $\R^{n}$. Without loss of generality, consider the case
$n=2$. The vector valued function $\vf(x)=\bigl(f_{1}(x),f_{2}(x)\bigr)^{\top}$
can be interpreted as the curve $C:I\ni x\mapsto(f_{1}(x),f_{2}(x),x)^{\top}\subset\R^{3}$.
Considering a partition $T$ of the interval $I$ into $N$ subintervals
$I_{i}=[x_{i-1},x_{i}]$, the linear interpolant between points $C(x_{i-1})$
and $C(x_{i})$ is the line joining them, i.e., $r_{i}:I_{i}\ni x\mapsto(\pip f_{1}(x),\pip f_{2}(x),x)^{\top}\subset\R^{3}$.
If, in each subinterval $I_{i}$, we measure the distance between
$C$ and the approximating line $r_{i}$ by the $\Lone$ distance
of their canonical projections on the first $n=2$ coordinate planes,
$\|\vf-\pip\vf\|_{\Lone(I_{i})}=\sum_{j=1}^{n}\|f_{j}-\pip f_{j}\|_{\Lone(I_{i})},$
where $\pip\vf(x)=(\pip f_{1}(x),\pip f_{2}(x))^{\top}$, then the
total distance between the curve $C$ and the polygonal line approximating
it (consisting of linear segments $\{r_{i}\}_{i=1}^{N}$) is
\begin{equation}
\|\vf-\pip\vf\|_{\Lone(I)}=\sum_{i=1}^{N}\|\vf-\pip\vf\|_{\Lone(I_{i})}.\label{eq:VectorValuedLOneDistance}
\end{equation}
In this framework, we may follow similar steps as those in Sections~\ref{sub:Error-Single-Interval}
through \ref{sub:Approximation-to-the-optimal-partition}, to conclude
that now $\|\vf''(x)\|_{1}$ plays the role of $|f''(x)|$. Hence~(\ref{eq:limit-optimal-vertex-density})
becomes $lkd(x)\propto\|\vf''(x)\|_{1}^{1/3}$ and this may be used
in the procedure~(\ref{eq:optimized-partition})-(\ref{eq:cumulative-knot-distribution})
to obtain the suboptimal partition and approximate upper bound formulas
(analogous to~(\ref{eq:ApproxL1Boundpip}) and~(\ref{eq:ApproxL1boundUniformPip}))
for~(\ref{eq:VectorValuedLOneDistance}). For large $N$, the $8/3$
factor between the errors of the vector valued linear interpolant
and the best $\Lone$ approximation still holds, i.e., $\|\vf-\hat{\vf}\|_{\Lone(I)}\approx\frac{3}{8}\|\vf-\pip\vf\|_{\Lone(I)}$
due to the independence of the optimizations in each coordinate plane
given a partition~$\part$.

\section{Computational complexity}

\label{sec:Complexity}The evaluation of $v(x)$ for any polygonal
function $v\in\vecspace$, such as $\fbest$ and $\pip\f$ previously
discussed, is very simple and consists of three steps: determination
of the subinterval $I_{i}=\left[x_{i-1},x_{i}\right)$ such that $x_{i-1}\leq x<x_{i}$,
computation of the fractional distance $\delta_{i}(x)=(x-x_{i-1})/(x_{i}-x_{i-1})$,
and interpolation of the resulting value $v(x)=\left(1-\delta_{i}(x)\right)v_{i-1}+\delta_{i}(x)v_{i}$.
Regardless of the specific polygonal function under consideration,
the computational cost of its evaluation is dominated by the first
step, which ultimately depends on whether or not the partition $\part$
is uniform. In the general case of $\part$ not being uniform, the
first step of the evaluation implies searching $\part$ for the correct
index $i$; since $\part$ is an ordered set, we can employ a binary
search to determine $i$, which means that the computational complexity
of the evaluation of $v(x)$ is $O(\log N)$ in the worst case. However,
in the particular case of $\part$ being uniform, the first and second
steps of the algorithm greatly simplify: $i\Leftarrow1+\left\lfloor N\left(x-x_{0}\right)/\left(x_{N}-x_{0}\right)\right\rfloor $
and $\delta_{i}(x)=1-i+N\left(x-x_{0}\right)/\left(x_{N}-x_{0}\right)$;
therefore, it suffices to store the endpoints $\{x_{0},x_{N}\}$ and,
most importantly, the computational complexity of the evaluation of
$v(x)$ becomes~$O(1)$. 

Consequently, approximations based on uniform partitions are expected
to perform better, in terms of execution time, than those based on
optimized partitions. However, if $x$ can be reasonably predicted
(e.g., due to it being the next sample of a well-characterized input
signal, such as in digital predistorters for power amplifiers~\cite{Cavers1999,Muhonen2000}),
other search algorithms with less mean computational complexity than
binary search could be used to benefit from the reduced memory requirement
of optimized partitions without incurring too great a computational
penalty.

The proposed algorithm is very simple to implement on either CPUs
or GPUs. However, the GPU case is specially relevant because its texture
filtering units are usually equipped with dedicated circuitry that
implements the interpolation step of the algorithm in hardware~\cite{doggett2012texture},
further accelerating evaluation.

\section{Experiments\label{sec:Experiments}}

\begin{figure}
\begin{centering}
\includegraphics[width=0.9\columnwidth]{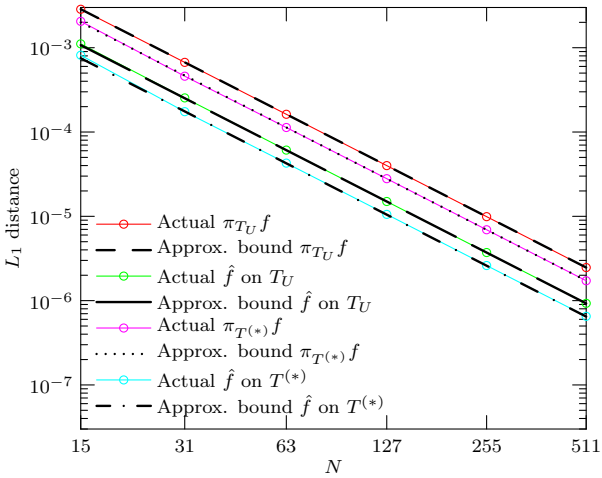}
\par\end{centering}

\caption{\label{fig:GraphicalErrors04}$\Lone$ distance for different polygonal
approximations to the Gaussian function $f(x)=\exp(-x^{2}/2)/\sqrt{2\pi}$
in the interval $x\in[0,4].$}
\end{figure}
To assess the performance of the developed linearization technique
we have selected a nonlinear function that is used in many applications
in every field of science, the Gaussian function $f(x)=\exp(-x^{2}/2)/\sqrt{2\pi}$.
We approximate it using a varying number of segments and a varying
domain.

Figure~\ref{fig:GraphicalErrors04} shows $\Lone$ distances between
the Gaussian function and the polygonal approximations described in
previous sections, in the interval $x\in[0,4]$. The best $\Lone$
polygonal approximation was computed as explained in Section~\ref{sub:Best-LOne-polygonal-computation},
using the Newton-Raphson iteration starting from the coefficients
of the best $\Ltwo$ approximation. The implementation relied in the
MATLAB function \texttt{\small{}lsqnonlin} with tridiagonal Hessian
and tolerances for stopping criteria: $10^{-16}$ in parameter space
and $10^{-10}$ in destination space. Convergence was fast, requiring
few iterations (typically less than 25 function evaluations) in a
process that is carried out once and off-line previous to the application
of the linearized function.

The computation of the optimized partition $\partOpt$ required independently
solving equation~(\ref{eq:optimized-partition}) for each of the
$N-1$ free knots of $\partOpt$. This solely relies in standard numerical
integration techniques, taking few seconds to complete, as opposed
to recursive partitioning techniques such as~\cite{Butler2011},
which take significantly longer.

The figure reports measured distances as well as the approximate upper
bounds to the distances~(\ref{eq:ApproxL1Boundpip}) and~(\ref{eq:ApproxL1boundUniformPip})
using the Riemann integral to approximate the sums. The measured $\Lone$
distances between $\pip f$ and $f$ using the uniform and optimized
partitions agree well with~(\ref{eq:ApproxL1boundUniformPip}) and~(\ref{eq:ApproxL1Boundpip}),
respectively, which have a $O(N^{-2})$ dependence rate that is also
applicable to the rest of the curves since they all have similar slopes.
The fit is good even for modest values of $N$ (e.g.,  $N=15$). Also,
the ratio between the distances corresponding to $\fbest$ and $\pip f$
is approximately the value $3/8$ originating from~(\ref{eq:approximate-best-L1-segment})
and~(\ref{eq:approximate-L1-interpolated-segment}). 

Equations~(\ref{eq:ApproxL1Boundpip}), (\ref{eq:ApproxL1boundUniformPip})
and/or the curves in Fig.~\ref{fig:GraphicalErrors04} can be used
to select the optimal value of $N$ solely based on distance considerations.
For example, in an application with a target $\Lone$ error tolerance
$\|\cdot\|_{\Lone(I)}\leq10^{-5}$, we obtain $N\ge\left\lceil 253.9\right\rceil =254$
for $\pipUnif f$ (Eq.~(\ref{eq:ApproxL1boundUniformPip})), $N\ge\left\lceil 212.2\right\rceil =213$
for $\pipOpt f$ (Eq.~(\ref{eq:ApproxL1Boundpip})), $N\ge\lceil253.9\sqrt{3/8}\rceil=\left\lceil 155.5\right\rceil =156$
for $\fbest$ on $\partUnif$ and $N\ge\lceil212.2\sqrt{3/8}\rceil=\left\lceil 129.9\right\rceil =130$
for $\fbest$ on~$\partOpt$. 

\begin{table}[t]
\caption{\label{tab:Mean-per-evaluation-execution}Mean per-evaluation execution
times (in picoseconds).\protect \\
$V_{\partUnif}$: polygonal functions defined on a uniform partition.\protect \\
$V_{\partOpt}$: polygonal functions defined on an optimized partition.}

\centering{}%
\begin{tabular}{|l|l|c|c|c|c|c|}
\hline 
\multicolumn{2}{|l|}{{\footnotesize{}Number of points ($N+1$)}} & {\footnotesize{}32} & {\footnotesize{}64} & {\footnotesize{}128} & {\footnotesize{}256} & {\footnotesize{}512}\tabularnewline
\hline 
\hline 
\multirow{3}{*}{{\footnotesize{}CPU}} & {\footnotesize{}Gaussian function } & \multicolumn{5}{c|}{{\footnotesize{}13710}}\tabularnewline
\cline{2-7} 
 & {\footnotesize{}Function in $V_{\partUnif}$} & \multicolumn{5}{c|}{{\footnotesize{}1750 - 1780}}\tabularnewline
\cline{2-7} 
 & {\footnotesize{}Function in $V_{\partOpt}$} & {\footnotesize{}8900} & {\footnotesize{}11230} & {\footnotesize{}13830} & {\footnotesize{}16910} & {\footnotesize{}20120}\tabularnewline
\hline 
\multirow{3}{*}{{\footnotesize{}GPU}} & {\footnotesize{}Gaussian function} & \multicolumn{5}{c|}{{\footnotesize{}14.2}}\tabularnewline
\cline{2-7} 
 & {\footnotesize{}Function in $V_{\partUnif}$} & \multicolumn{5}{c|}{{\footnotesize{}7.8 - 7.9}}\tabularnewline
\cline{2-7} 
 & {\footnotesize{}Function in $V_{\partOpt}$} & {\footnotesize{}122.7} & {\footnotesize{}142.9} & {\footnotesize{}163.2} & {\footnotesize{}188.3} & {\footnotesize{}211.1}\tabularnewline
\hline 
\end{tabular}
\end{table}
Table~\ref{tab:Mean-per-evaluation-execution} shows mean processing
times per evaluation both on a CPU (sequentially, one core only) and
on a GPU. All execution time measurements have been taken in the same
computer, equipped with an Intel Core i7-2600K processor, 16 GiB RAM
and an NVIDIA GTX 580 GPU. We compare the fastest option~\cite{NVIDIACorporation2012}
for implementing the Gaussian function using its definition against
its approximation using the proposed algorithm. Note that the processing
time of any polygonal function $v\in V_{T}$ solely depends on $T$,
as shown in Section~\ref{sec:Complexity}; as expected from the analysis,
execution times are constant in the case of a uniform partition and
grow logarithmically with $N$ in the case of an optimized partition.
The proposed strategy, using a uniform partition, solidly outperforms
conventional evaluation of the nonlinear function.

\begin{figure}[t]
\begin{centering}
\includegraphics[width=0.9\columnwidth]{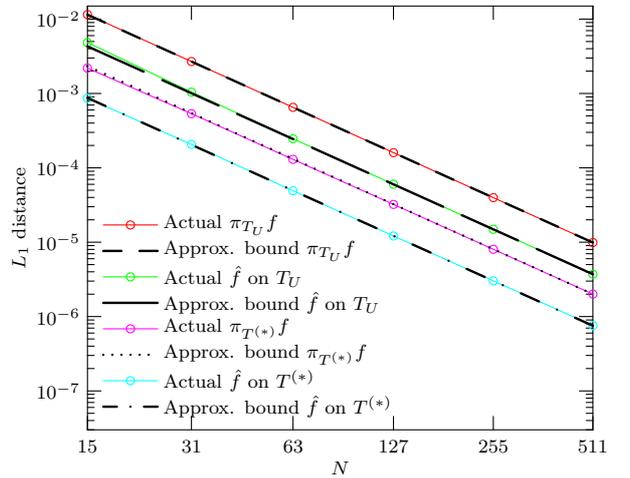}
\par\end{centering}

\caption{\label{fig:GraphicalErrors08}$\Lone$ distance for different polygonal
approximations to the Gaussian function $f(x)=\exp(-x^{2}/2)/\sqrt{2\pi}$
in the interval $x\in[0,8].$}
\end{figure}
\begin{figure}[t]
\begin{centering}
\includegraphics[width=0.8\columnwidth]{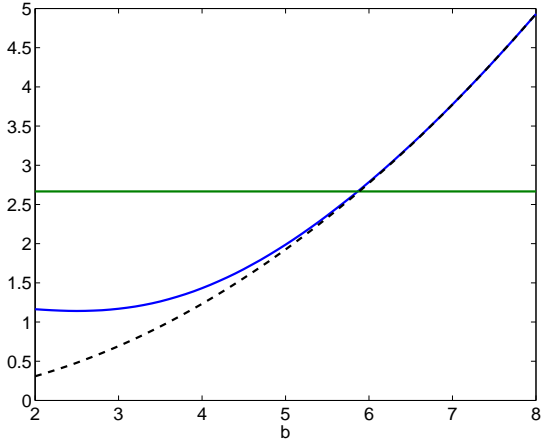}
\par\end{centering}

\caption{\label{fig:Gain3to8}Blue: gain obtained by optimizing a partition:
(\ref{eq:ApproxL1boundUniformPip})/(\ref{eq:ApproxL1Boundpip}) for
the Gaussian function $f(x)=\exp(-x^{2}/2)/\sqrt{2\pi}$, as a function
of $b$ in the interval $I=[0,b]$. Green: gain obtained by using
the best $\Lone$ approximation instead of the interpolant.}
\end{figure}
The approximation errors in $\Lone$ distance over the interval $x\in[0,8]$
were also measured. These measurements and the corresponding approximate
upper bounds are reported in Fig.~\ref{fig:GraphicalErrors08}. Observe
that, in this case, the curve $\|\f-\fbest_{\partUnif}\|_{\Lone([0,8])}$
is above $\|\f-\pipQuasiOpt\f\|_{\Lone([0,8])}$, whereas in the interval
$[0,4]$ the relation is the opposite (Fig.~\ref{fig:GraphicalErrors04}).
This issue is easily explained by our previous error analysis: the
gap between the approximation errors of the interpolant and the best
$\Lone$ approximation is constant ((\ref{eq:approximate-L1-interpolated-segment})/(\ref{eq:approximate-best-L1-segment})$\approx8/3$),
whereas the gain obtained by optimizing a partition ($\|\f-\pipUnif f\|_{\Lone(I)}/\|\f-\pipOpt f\|_{\Lone(I)}$
or $\|\f-\fbest_{\partUnif}\|_{\Lone(I)}/\|\f-\fbest_{\partOpt}\|_{\Lone(I)}$)
depends on the approximation domain $I$. Fig.~\ref{fig:Gain3to8}
shows both gains as functions of $b$ in the interval $I=[a,b]$,
with $a=0$. The horizontal line at $8/3$ corresponds to the gain
obtained by using the best $\Lone$ approximation instead of the interpolant,
regardless of the partition. The blue solid line shows the gain obtained
by optimizing a partition, (\ref{eq:ApproxL1boundUniformPip})/(\ref{eq:ApproxL1Boundpip}).
As $b$ increases, it behaves asymptotically as the parabola (\ref{eq:ApproxL1boundUniformPip})/(\ref{eq:ApproxL1Boundpip})~$\approx0.077b^{2}$
(dashed line), which can readily be seen by taking the limit $\lim_{b\to\infty}\|\fsec\|_{\Lone(I)}/\bigl(\int_{a}^{b}\left|\fsec(t)\right|^{1/3}\mathrm{d}t\bigr)^{3}\approx0.077$.
As $b$ increases, most of the gain is due to the approximation of
the tail of the Gaussian by few and large linear segments, which leaves
most of the budget of control points to better approximate the central
part of the Gaussian. The point at which the gain (blue line) meets
the horizontal line at $8/3$ indicates the value of $b$ where the
$\|\f-\fbest_{\partUnif}\|_{\Lone([0,b])}$ and $\|\f-\pipQuasiOpt\f\|_{\Lone([0,b])}$
curves swap positions.

\subsection*{Chirp function}

\begin{figure}
\centering{}\includegraphics[width=1\columnwidth]{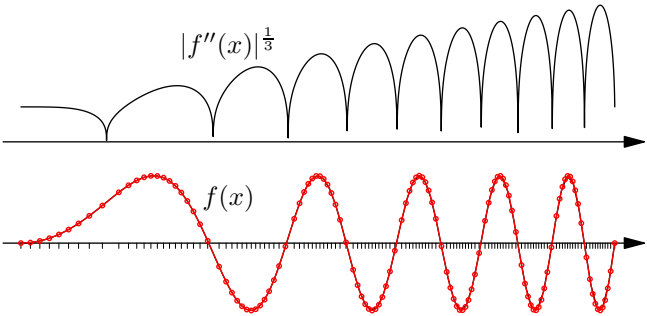}\caption{\textbf{\label{fig:Chirp}}Suboptimal partition for the linear chirp
function $f(x)=\sin(10\pi x^{2})$ in the interval $x\in[0,1]$. Top:
local knot density ($lkd$) corresponding to $f$. Bottom: polygonal
interpolant $\pipOpt f$ with $N=127$ ($128$ knots) overlaid on
function $f$. }
\end{figure}

\begin{figure}
\begin{centering}
\includegraphics[width=0.9\columnwidth]{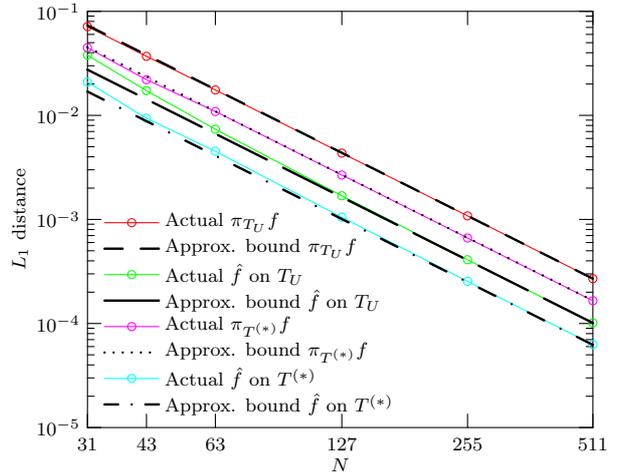}
\par\end{centering}

\caption{\label{fig:GraphicalErrorsChirp}$\Lone$ distance for different polygonal
approximations to the chirp function $f(x)=\sin(10\pi x^{2})$ in
the interval $x\in[0,1]$.}
\end{figure}

The performance of the linearization technique has also been tested
on a more challenging case: the chirp function $f(x)=\sin(10\pi x^{2})$,
which combines both nearly flat and highly oscillatory parts (see
Fig.~\ref{fig:Chirp}). Figure~\ref{fig:GraphicalErrorsChirp} reports
the $\Lone$ distances between the chirp function and the polygonal
approximations described in previous sections, in the interval $x\in[0,1]$.
For $N\geq63$ the measured errors agree well with the predicted approximate
error bounds, whereas for $N<63$ the measured errors differ from
the predicted ones (specifically in the optimized partition) because
in these cases the number of linear segments is not enough to properly
represent the high frequency oscillations. 

A sample optimized partition and the corresponding polygonal interpolant
$\pipOpt f$ is also represented in Fig.~\ref{fig:Chirp}. The knots
of the partition are distributed according to the local knot density
($lkd$) (see Fig.~\ref{fig:Chirp}, Top), whose envelope grows according
to $x^{2/3}$, and this trend is reflected in the accumulation of
knots in the regions of high oscillations (excluding the places around
the zeros of the $lkd$). 

The evaluation times coincide with those of the Gaussian function
(Table~\ref{tab:Mean-per-evaluation-execution}) because the processing
time of the polygonal approximation does not depend on the function
values.

\section{Conclusions\label{sec:Conclusions}}

We have developed a practical method to linearize and numerically
evaluate arbitrary continuous real-valued functions in a given interval
using simpler polygonal functions and measuring errors according to
the $\Lone$ distance. As a by-product, our technique allows fast
(e.g., real-time) implementation of computationally expensive applications
that use such mathematical functions. 

To this end, we analyzed the polygonal approximations given by the
linear interpolant and the least-first-power or best $\Lone$ approximation
of a function. A detailed error analysis in the $\Lone$ distance
was carried out seeking a nearly optimal design of both approximations
for a given budget of subintervals~$N$. In the practical asymptotic
case of large $N$, we used error equalization to achieve a suboptimal
design (partition~$\part$) and derive a tight bound on the approximation
error for the linear interpolant, showing $O(N^{-2})$ dependence
rate that was confirmed experimentally. The best $\Lone$ approximation
improves upon the results of the linear interpolant by a rough factor
of 8/3. 

Combining both quality and computational cost criteria, we conclude
from this investigation that, from an engineering standpoint, using
the best $\Lone$ polygonal approximation in uniform partitions is
an excellent choice: it is simple, fast and its error performance
is very close to the limit defined by optimal partitions. Possible
paths to explore related to our technique are, among others, extension
to multidimensional functions and the incorporation of constraints
in the linearization process (e.g., so that the best $\Lone$ polygonal
model also satisfies positivity or a target bounded range). 

%\bibliographystyle{IEEEtran}
%\bibliography{gaussianasL1}

\end{document}